# A remark on extension of log-pluricanononical forms for nondegenerate hypersurfaces

Achim Hennings [♣]

**Abstract**: We determine, complementing [M], the log-pluricanonical forms on a nondegenerate hypersurface. This description shows that they are extendable to 1-parameter deformations. For an equisingular deformation we thus obtain a simultaneous log-canonical model.

Let $f \in \mathbb{C}[x_1, \ldots, x_n] = P$ be a polynomial satisfying $f(\ldots, 0, x_i, 0, \ldots) \neq 0$ for all variables. Let $\Delta = \Gamma_+(f)$ be its Newton polyhedron, $\Gamma$ the compact part of the boundary of $\Delta$, and $\Sigma_0$ the fan dual to $\Delta$. Let $V = (f) \subseteq \mathbb{C}^n$ be the divisor of $f$. We wish to determine all pluricanonical forms on (the smooth part of) $V$, which have at most logarithmic poles on a resolution of $V$ with normal crossings of the exceptional set. This condition is independent of the resolution. We show that, if $f$ is nondegenerate, these log-pluricanonical forms can be obtained as residues from a toric embedded resolution and have a simple description in terms of $\Delta$. In [M] this and more is claimed, but the proofs appear incomplete. Subsequently, we use this result to prove that log-pluricanonical forms admit extensions to 1-parameter deformations of $V$. If an equisingular deformation (in the sense of simultaneous resolution) is given, one thereby obtains a simultaneous log-canonical model of the deformation.

## 1 Log-pluricanonical forms

We show that these forms are the residues of the corresponding forms on a toric embedded resolution. Therefore, the filtration on $P/fP$ by log-canonical forms can be immediately described by integer multiples of $\Delta$. (Cf. [M].)

We consider the boundary hyperplanes $H_i = \{r \in \mathbb{R}^n | r_i = 0\}$ and the intersections $\Delta \cap H_i$. We may take, e.g., $i = 1$. Let $\delta_j$ denote the compact facets $\delta$ of $\Delta$ which meet $H_1$. They correspond to the primitive normal vectors $v_j$ that are together with $e_1$ among the edges of a cone in $\Sigma_0$. We denote by $\lambda_j$ the linear form $\langle v_j, . \rangle$ on $\mathbb{R}^n$ (considered as the vector space containing $\Delta$). We form the convex polyhedron

$$\Delta_1 = \{r \in \mathbb{R}^n_+ | \lambda_j(r) \geq \min \lambda_j(\Delta) \; \forall j\}$$

with $\Gamma_1$ the compact part of the boundary. Since $\Delta$ is an integral polyhedron,

$$\Delta_1 \cap [0,1] \times \mathbb{R}^{n-1}_+ = \Delta \cap [0,1] \times \mathbb{R}^{n-1}_+.$$

Each compact facet $\delta$ of $\Delta_1$ contains a facet $\delta_j$ and therefore meets $H_1$. Namely, $\delta$ is represented as $\Delta_1 \cap \lambda_j^{-1}(\min \lambda_j(\Delta))$ for some $j$. The polyhedra $\Delta$ and $\Delta_1$ define rational decreasing filtrations on $P$:

$$F^a = L(a\Delta), F_1^a = L(a\Delta_1).$$

Here we use $L(A)$ to denote the vector space of (Laurent-) polynomials with support (i.e. exponents) in the set $A$, and $F^a = F_1^a = P$ for $a \leq 0$ by definition. We have

$$F^a \subseteq F_1^a \text{ and } F_1^a F_1^b \subseteq F_1^{a+b} \text{ (and the same for } F\text{)}.$$

[♣]Universität Siegen, Fakultät IV, Hölderlinstraße 3, D-57068 Siegen

**Lemma 1:** The multiplication map ($a \geq 0$, $k \in \mathbb{N}$)

$$P/\left(F_1^a + (x_1^k)\right) \xrightarrow{f} P/\left(F_1^{a+1} + (x_1^k)\right)$$

is injective.

Proof: It is sufficient, if $Gr_{F_1}^a(P/x_1^k P) \xrightarrow{f} Gr_{F_1}^{a+1}(P/x_1^k P)$, $a \geq 0$, is injective. Let $g \neq 0$ be a polynomial of $x_1$−degree $\leq k-1$ with support in $a\Gamma_1$. Let $\delta$ be a compact facet of $\Delta_1$ with $\tilde{g} = g_{a\delta} \neq 0$ (the part of $g$ with support in $a\delta$). Let $\delta' = \delta \cap \Delta$ be the corresponding boundary facet of $\Delta$ and let $f_{\delta'}$, resp. $f_{\delta' \cap H_1}$ the part of $f$ belonging to $\delta'$ resp. $\delta' \cap H_1$. As $f_{\delta' \cap H_1} \neq 0$, we get $(f_{\delta'} \tilde{g} \mod x_1^k P) \neq 0$, by considering the lowest term with respect to $x_1$ in $\tilde{g}$. This product is also the $(a+1)\delta$−part of the product $[fg]$ in $Gr_{F^1}^{a+1}(P/x_1^k P)$, and therefore $[fg] \neq 0$.

**Lemma 2:** The filtration $F$ has the property: $gf \in F^a + x_1^k P$ implies $g \in F^{a-1} + x_1^k P$, provided $k \leq a$.

Proof: For polynomials $h$ of $x_1$−degree $\leq a$ we have $h \in F_1^a \Leftrightarrow h \in F^a$ because $a\Delta_1 \cap [0,a] \times \mathbb{R}_+^{n-1} = a\Delta \cap [0,a] \times \mathbb{R}_+^{n-1}$. So we get:

$$gf \in F^a + x_1^k P \Longrightarrow gf \in F_1^a + x_1^k P$$

$$\Longrightarrow g \in F_1^{a-1} + x_1^k P \quad \text{(Lemma 1)}$$

$$\Longrightarrow g \in F^{a-1} + x_1^k P \quad \text{(taking } g \text{ of degree} \leq k-1 \text{ in } x_1\text{)}.$$

**Lemma 3:** $gf \in F^a + (x_1 \ldots x_n)^k$ implies $g \in F^{a-1} + (x_1 \ldots x_n)^k$, provided $k \leq a$.

Proof: By Lemma 2 we have $g \in F^{a-1} + x_i^k P$ $\forall i$, thus $g \in F^{a-1} + (x_1 \ldots x_n)^k$.

From now on, we assume that $f$ is nondegenerate (in the sense of [K]). For simplicity, we may assume that $V^* = V - \{0\}$ is non-singular and the same holds for the projective closure (by adding a suitable homogeneous polynomial of high degree). To avoid special cases we also assume $n \geq 3$.

Let $\pi: X = X(\Sigma) \to \mathbb{C}^n$ be a resolution of $f$ associated to a regular subdivision $\Sigma$ of $\Sigma_0$. Let $D = \bigcup_j D_j$ be the exceptional set, $D_i'$ the strict transform of $(x_i)$, $D' = \bigcup_i D_i'$, and $(\pi^* f) = M + C$ the decomposition into strict transform and exceptional part. We put $E = D \cap M$ and $E' = D' \cap M$.

The Poincaré residue defines a map

$$\omega_X(D+M)^m \to \omega_M(E)^m.$$

In terms of divisors, this is the map

$$\mathcal{O}_X(m(K+D+M)) \to \mathcal{O}_M(m(K+D+M)).$$

Here $K + D + D' \sim 0$, so $K + D \sim -D'$.

**Lemma 4:**

$$H^0\left(\mathcal{O}_X(m(K+D+M))\right) \to H^0\left(\mathcal{O}_M(m(K+D+M))\right)$$

is surjective ($m \geq 1$).

Proof: We have (since $C, D'$ are coprime and $\mathcal{O}_X(-(m-1)C)$ is acyclic)

$$H^0\left(\mathcal{O}_X(m(K+D+M))\right) \cong H^0\left(\mathcal{O}_X(m(-D'-C))\right)$$
$$= H^0(\mathcal{O}_X(-mC)) \cap H^0(\mathcal{O}_X(-mD'))$$
$$= L(m\Delta) \cap (x_1 \ldots x_n)^m \subseteq P,$$

$$H^0\left(\mathcal{O}_M(m(K+D+M))\right) \cong H^0\left(\mathcal{O}_M(m(-D'-C))\right)$$
$$= H^0(\mathcal{O}_M(-mC)) \cap H^0(\mathcal{O}_M(-mD'))$$
$$= F^m(P/fP) \cap (x_1 \ldots x_n)^m P/fP \subseteq P/fP.$$

Let $h \in F^m$ be a representative for some element of the last intersection.

$$\Rightarrow \exists\, g \in P\colon h - gf \in (x_1 \ldots x_n)^m$$
$$\Rightarrow gf \in F^m + (x_1 \ldots x_n)^m$$
$$\Rightarrow g \in F^{m-1} + (x_1 \ldots x_n)^m \quad \text{(Lemma 3)}.$$

Write $g = g_1 + g_2$ with $g_1 \in F^{m-1}$, $g_2 \in (x_1 \ldots x_n)^m$. Then $h - g_1 f \in F^m$ and $h - g_1 f = (h - gf) + g_2 f \in (x_1 \ldots x_n)^m$. Thus $h - g_1 f \in F^m \cap (x_1 \ldots x_n)^m$ is a new representative in $H^0\left(\mathcal{O}_X(m(-D'-C))\right)$.

## 2 Extension of log-pluricanonical forms

Let $\mathcal{V} \subseteq \mathbb{C}^n \times T$ ($T \subseteq \mathbb{C}$ a disc about 0) be a deformation of $V$, defined by a polynomial $F(x,t)$ with $F(x,0) = f(x)$. We may assume that $\mathcal{V} \to T$ is smooth at infinity by taking $F$ of lower degree. Let $\mathcal{V}^*$ be the smooth locus of $\mathcal{V}$. We show that the log-pluricanonical forms on $V^*$ can be extended to such forms on $\mathcal{V}^*$.

Let $\pi\colon X \to \mathbb{C}^n$ be the toric resolution from section 1. Let $\alpha \in H^0(V^*, \omega_V^m)$ be an $m$-fold log-canonical form. By Lemma 4, $\alpha$ has a representation as an $m$-fold Poincaré residue

$$\operatorname{Res}_{V^*}^{[m]} \frac{g(dx/x)^m}{f^m} = \operatorname{Res}_{V^*}^{[m]} \frac{h(dx)^m}{f^m}$$

with $g \in L(m\Delta) \cap (x_1 \ldots x_n)^m$, $h = (x_1 \ldots x_n)^{-m} g$. Here

$$\beta = \pi^* \frac{h(dx)^m}{f^m} \in H^0(X, \omega_X(D+M)^m)$$

and $\operatorname{Res}_M^{[m]} \beta \in H^0(M, \omega_M(E)^m)$. Since $\beta$ has poles of order at most $m$ along $D_j$ and the order of $f$ is positive, for any $\lambda > 0$ the form

$$|\pi^* f|^{2\lambda} (\beta \wedge \bar{\beta})^{1/m} = \text{const} \cdot |\pi^* f|^{2(\lambda-1)} |\pi^* h|^{2/m} \pi^*(dx \wedge d\bar{x})$$

is locally integrable on $X$. For a ball $B \subseteq \mathbb{C}^n$ around 0 it follows then that

$$\int_B |f|^{2(\lambda-1)} |h|^{2/m} |dx \wedge d\bar{x}| < \infty, \quad \lambda > 0$$

(where $|dx \wedge d\bar{x}|$ denotes the positive measure of the $2n$-form).

We apply the extension theorem of B. Berndtsson und M. Paun [P]. Accordingly, $h$ can be extended to a holomorphic function $H(x,t)$ such that

$$\int_{B\times T}|F|^{2(\lambda-1)}|H|^{\frac{2}{m}}|d(x,t)\wedge d(\bar{x},\bar{t})|<\infty,$$

where $H$ depends on $\lambda > 0$.[1]

Now let $\rho: Y \to B \times T$ be a resolution of $F$. This means that the exceptional set $A = \bigcup_k A_k$ and the strict transform $N$ of $\mathcal{V}$ form a simple normal crossings divisor. We may write $(\rho^*F) = N + \sum_k a_k A_k$, $(\rho^*d(x,t)) = \sum_k b_k A_k$. We choose $0 < \lambda < \frac{1}{ma_k}$ $\forall k$. Let $H$ be an extension of $h$ as above. Then $F^{\lambda-1}H^{\frac{1}{m}}d(x,t)$ has pole order $< 1$ along $A_k$ $\forall k$, i.e.

$$(\lambda-1)a_k + \frac{1}{m}v_{A_k}(\rho^*H) + b_k > -1$$

$$\Leftrightarrow \frac{1}{m}v_{A_k}(\rho^*H) + b_k - a_k > -1 - \lambda a_k.$$

The left term is in $\frac{1}{m}\mathbb{Z}$, and $\lambda a_k \in ]0, \frac{1}{m}[$. Hence

$$\frac{1}{m}v_{A_k}(\rho^*H) + b_k - a_k \geq -1, \text{ i.e. } v_{A_k}\left(\rho^*\frac{Hd(x,t)^m}{F^m}\right) \geq -m.$$

Thus

$$Res_N^{[m]}\rho^*\frac{Hd(x,t)^m}{F^m} \in H^0(N, \omega_N(A\cap N)^m)$$

is an extension of $\alpha$ to $\mathcal{V}^*$ with only logarithmic poles.

**Theorem:** Any log-pluricanonical form on $(V, 0)$ can be extended to $(\mathcal{V}, 0)$.

Remark: According to S. Ishii [I1], on a normal Gorenstein singularity $V$, the $L^{2/m}$ −integrable $m$ −fold canonical forms are always extendable. These are the $m$ −fold canonical forms on a good resolution with pole order $\leq m-1$ (as opposed to $\leq m$ here).

**Corollary:** Let $\rho: \mathcal{X} \to \mathcal{V}$ be a projective weak simultaneous resolution, i.e. $\mathcal{X}_t \to \mathcal{V}_t$ is a resolution and the exceptional set $\mathcal{E}$ is locally trivial over $T$. Then $\mathcal{V}$ has a simultaneous log-canonical model.

Proof: We can assume that $\mathcal{E}$ has normal crossings over $T$. The log-canonical model of $\mathcal{V}$ resp. $V$ exists (by [OX]) and is given by

$$\mathcal{Z} = Proj \bigoplus_{m\geq 0} \rho_*\mathcal{O}(m(K_\mathcal{X} + \mathcal{E})) \text{ resp.}$$

$$Z = Proj \bigoplus_{m\geq 0} \rho_{0*}\mathcal{O}\left(m(K_{\mathcal{X}_0} + \mathcal{E}_0)\right).$$

By the theorem the graded algebra of $\mathcal{Z}$ is flat over $T$ and it is mapped surjectively onto the graded algebra of $Z$ with kernel generated by $t$. Thus, $Z$ is the special fibre of the flat map $\mathcal{Z} \to T$.

---

[1] Note that $\log|F|$ is plurisubharmonic and the theorem is applied for $\lambda < 1$